# Both real and imaginary parts of the function $F(s)=\varsigma(s)\Gamma(s/2)\pi^{-s/2}=\frac{\xi(s)}{s(s-1)}$, whose zeroes exactly coincide with the non-trivial zeroes of the Riemann $\varsigma$-function, have infinitely many zeroes for any value of *Re*s.

The following theorem is proven: Both real and imaginary parts of the function $F(s)$ defined as $F(s)=\varsigma(s)\Gamma(\frac{1}{2}s)\pi^{-s/2}=\frac{\xi(s)}{s(s-1)}$, and whose zeroes exactly coincide with the non-trivial zeroes of the Riemann $\varsigma$-function, have infinitely many zeroes for any value of Re$s$.

S. K. Sekatskii, Laboratoire de Physique de la Matière Vivante, IPSB, BSP 408, Ecole Polytechnique Fédérale de Lausanne, CH1015 Lausanne-Dorigny, Switzerland.
E-mail : serguei.sekatski@epfl.ch

This is evident that zeroes of the Riemann $\xi$-function $\xi(s) = \varsigma(s)\Gamma(\frac{1}{2}s)\pi^{-s/2}s(s-1)$ coincide with the non-trivial zeroes of the Riemann $\varsigma$-function (see e.g. [1, 2]). Of course, the same is valid for the function $F(s) = \varsigma(s)\Gamma(s/2)\pi^{-s/2} = \frac{\xi(s)}{s(s-1)}$ which, contrary to the Riemann $\xi$-function is not a whole function but has simple poles at $z=0$ and $z=1$ with the residues equal respectively to -1 and 1. Nevertheless, we believe that introduction of this function is still useful to shorten the notation of the paper. The aim of the present Note is the proof of the following

THEOREM. Both real and imaginary parts of the function $F(s)$ defined as $F(s) = \varsigma(s)\Gamma(\frac{1}{2}s)\pi^{-s/2} = \frac{\xi(s)}{s(s-1)}$, and whose zeroes exactly coincide with the non-trivial zeroes of the Riemann $\varsigma$-function, have infinitely many zeroes for any value of Re$s$.

Below we will use also the notation $F(s) \equiv F_\alpha(t)$, where $s = \sigma + it = \frac{1}{2} + \alpha + it$ with $\sigma, \alpha, t$ real. In terms of the function $F_\alpha(t)$ the Riemann hypothesis reads simply that there are no zeroes for $\alpha \neq 0$. Note also that $F_0(t)$ is real; of course $F_\alpha(t) = F^*_\alpha(-t)$ and the relation $F_\alpha(t) = F_{-\alpha}(-t)$ holds. The method of proof of our Theorem is a slight generalization of Hardy paper [3] where he showed that there are infinitely many Riemann function zeroes lying on the line $s = 1/2 + it$. Evidently, this is a particular case concerning $\operatorname{Re} F(s)$ for $s=1/2$ of our Theorem supplemented with the trivial $\operatorname{Im} F(s) = 0$ for $s=1/2$.

*Proof.* We start from $\dfrac{1}{4\pi i}\int_{k-i\infty}^{k+i\infty}\Gamma(u/2)(\pi x)^{-u/2}\zeta(u)du = w(x)$, where $w(x)=\sum_{n=1}^{\infty}\exp(-\pi n^2 x)$ and which is valid when $k>1$ and $\operatorname{Re} x>0$ [1-3] and introduce a contour integral $\int_C \Gamma(u/2)(\pi x)^{-u/2}\zeta(u)du$ round the rectangular contour C composed by straight lines connecting the points $1/2+\alpha-iT$, $1/2+\alpha+iT$, $k+iT$, $k-iT$ with $T\to+\infty$. We take $k>1$ and, for a moment, $-1/2<\alpha<1/2$. Of course, the integrand is nothing else than $F(u)x^{-u/2}$ hence we are speaking about the integral $\int_C F(u)x^{-u/2}du$. Inside the contour the integrand has a simple pole at $u=1$ with a residue equal to $x^{-1/2}$, on its left border we have an integral $-i\int_{-\infty}^{\infty}x^{-1/4-\alpha/2-it/2}F_\alpha(t)dt$ while on the right border the value of integral is equal to $4\pi i w(x)$, see above.

The disappearance of integrals taken over horizontal lines for $T\to+\infty$ is evident (again, similar to Hardy's case [3]) and thus Cauchy theorem gives

$$-i\int_{-\infty}^{\infty}x^{-1/4-\alpha/2-it/2}F_\alpha(t)dt + 4\pi i w(x) = 2\pi i x^{-1/2},$$ that is

$\dfrac{1}{2\pi}\int_{-\infty}^{\infty}x^{-1/4-\alpha/2-it/2}F_\alpha(t)dt = x^{-1/2} - 2w(x)$. Now taking $x=\exp(ib)$, $b$ real positive such that $\cos(b)>0$, one gets

$$\dfrac{1}{2\pi}\int_{-\infty}^{\infty}\exp(bt/2)F_\alpha(t)dt = e^{ib(1/4+\alpha/2)}(e^{-ib/2}-2w(e^{ib})) = e^{ib(-1/4+\alpha/2)}-2e^{ib(1/4+\alpha/2)}w(e^{ib}) =$$
$$e^{ib(-1/4+\alpha/2)} - e^{ib(1/4+\alpha/2)}(1+2w(e^{ib})) + e^{ib(1/4+\alpha/2)}$$

Taking $x=\exp(-ib)$ he or she gets

$$\dfrac{1}{2\pi}\int_{-\infty}^{\infty}\exp(-bt/2)F_\alpha(t)dt = e^{-ib(-1/4+\alpha/2)} - e^{-ib(1/4+\alpha/2)}(1+2w(e^{-ib})) + e^{-ib(1/4+\alpha/2)} =$$
$$e^{-ib(-1/4+\alpha/2)} - e^{-ib(-1/4+\alpha/2)}(1+2w(e^{ib})) + e^{-ib(1/4+\alpha/2)}$$

(Relation $1 + 2w(x) = \frac{1}{\sqrt{x}}(1 + 2w(1/x))$ ([2], P. 273), that is $1 + 2w(e^{-ib}) = e^{ib/2}(1 + 2w(e^{ib}))$, is used to obtain the last equality). By summing these two we have

$$\frac{1}{\pi}\int_{-\infty}^{\infty} \cosh(bt/2) F_\alpha(t) dt = 2\cos(b/4 - \alpha b/2) + 2\cos(b/4 + b\alpha/2) - \quad (8)$$
$$(e^{-ib(-1/4+\alpha/2)} + e^{ib(1/4+\alpha/2)})(1 + 2w(e^{ib}))$$

and by finding the difference we have

$$\frac{1}{\pi}\int_{-\infty}^{\infty} \sinh(bt/2) F_\alpha(t) dt = -2i\sin(b/4 - \alpha b/2) + 2i\sin(b/4 + b\alpha/2) - \quad (9)$$
$$(-e^{-ib(-1/4+\alpha/2)} + e^{ib(1/4+\alpha/2)})(1 + 2w(e^{ib}))$$

Evidently, (8) gives a correct Hardy's limit for $\alpha = 0$ [3]:

$$\frac{1}{\pi}\int_0^{\infty} \cosh(bt/2) F_0(t) dt = 2\cos(b/4) - e^{ib/4}(1 + 2w(e^{ib})).$$

Functions $\cosh(bt/2)$, $\sinh(bt/2)$ are real and thus we can write separating real and imaginary parts of (8), (9) and taking into account the odd/even nature of the functions involved:

$$\frac{1}{\pi}\int_0^{\infty} \cosh(bt/2)\,\mathrm{Re}\,F_\alpha(t) dt = \quad (10)$$
$$\cos(b/4 - \alpha b/2) + \cos(b/4 + \alpha b/2) - \frac{1}{2}\mathrm{Re}((e^{ib(1/4+\alpha/2)} + e^{-ib(-1/4+\alpha/2)})(1 + 2w(e^{ib})))$$

$$\frac{1}{\pi}\int_0^{\infty} \sinh(bt/2)\,\mathrm{Im}\,F_\alpha(t) dt = \quad (11)$$
$$-\sin(b/4 + \alpha b/2) + \sin(b/4 - \alpha b/2) - \frac{1}{2}\mathrm{Im}((e^{ib(1/4+\alpha/2)} - e^{-ib(-1/4+\alpha/2)})(1 + 2w(e^{ib})))$$

For further analysis we use exactly the same approach as Hardy did differentiating (8) and (9) *2p* times with respect to *b* and taking the limit $b \to \pi/2$. An asymptotic of the function $\Xi_\alpha(t)$ in any sector, $\Xi_\alpha(t) = O(t^A \exp(-\pi t/4))$ [1], enables to differentiate these equations with respect to *b* infinitely many times provided that $b < \pi/2$. The most important observation here is that the function $1 + 2w(e^{ib})$ and all its derivatives $\frac{d^p}{db^p}(1 + 2w(e^{ib}))$ tend to zero in this limit as this is known from the theory of elliptic functions and was used by Hardy; for short and clear proof see p. 215 of Titchmarsh and Heath-Brown book [1].

Thus we have in the limit $b \to \pi/2$:

$$\frac{1}{\pi 2^{2p}} \int_0^\infty t^{2p} \cosh(\pi t/4) \operatorname{Re} F_\alpha(t) dt = \\ (-1)^p \frac{1}{4^{2p}} [(1-2\alpha)^{2p} \cos\left(\frac{\pi}{8}(1-2\alpha)\right) + (1+2\alpha)^{2p} \cos\left(\frac{\pi}{8}(1+2\alpha)\right)] \quad (12)$$

and

$$\frac{1}{\pi 2^{2p}} \int_0^\infty t^{2p} \sinh(\pi t/4) \operatorname{Im} F_\alpha(t) dt = \\ (-1)^p \frac{1}{4^{2p}} [(1-2\alpha)^{2p} \sin\left(\frac{\pi}{8}(1-2\alpha)\right) - (1+2\alpha)^{2p} \sin\left(\frac{\pi}{8}(1+2\alpha)\right)] \quad (13)$$

From these expressions it immediately follows that their r.h.s change sign any time when *p* changes to *p+1*. (For (12) this is trivial because all factors in the square brackets are positive; for (13) it suffices to say that if $\alpha > 0$ than $\sin\left(\frac{\pi}{8}(1+2\alpha)\right) > \sin\left(\frac{\pi}{8}(1-2\alpha)\right)$ and $(1+2\alpha)^{2p} > (1-2\alpha)^{2p}$ thus the term in the square brackets is always negative; similarly for $\alpha < 0$). Hence we can use a variant of the Fejér theorem that the number of changes of the sign of a continuous function *f(x)* in a certain interval *(0, c)* is no less than the number of changes in sign of the sequence $f(0), \int_0^c f(t)dt, ..., \int_0^c f(t)t^n dt, ...$, cf. [1] P. 258, to prove our theorem for the case $-1/2 < \alpha < 1/2$.

The cases $\alpha < -1/2$ and $\alpha > 1/2$ are easily considered along the same lines. (Actually they are a bit easier). If, for example, $\alpha < -1/2$, in the interior of our contour we have additionally a simple pole at *u=0* with a residue equal to *-1*. Thus Cauchy theorem gives $\frac{1}{2\pi}\int_{-\infty}^{\infty}\exp(bt/2)F_\alpha(t)dt = e^{ib(1/4+\alpha/2)}(e^{-ib/2}-(1+2w(e^{ib})))$ from which repeating the same as above and taking the limit $b \to \pi/2$ we have after differentiation *2p* times:

$$\frac{1}{\pi 2^{2p}}\int_0^\infty t^{2p}\cosh(\pi t/4)\operatorname{Re} F_\alpha(t)dt = (-1)^p \frac{1}{4^{2p}}(1-2\alpha)^{2p}\cos\left(\frac{\pi}{8}(1-2\alpha)\right)$$

$$\frac{1}{\pi 2^{2p}}\int_0^\infty t^{2p}\sinh(\pi t/4)\operatorname{Im} F_\alpha(t)dt = (-1)^p \frac{1}{4^{2p}}(1-2\alpha)^{2p}\sin\left(\frac{\pi}{8}(1-2\alpha)\right).$$

Thus our theorem again follows apart from the cases where $\cos\left(\frac{\pi}{8}(1-2\alpha)\right)$ or $\sin\left(\frac{\pi}{8}(1-2\alpha)\right)$ are equal to zero. But if this is the case, we can differentiate (14) or (15) not *2p* but *2p+1* (odd number) times - cosine then changes to sine and vice versa, and the theorem follows again.

The case $\alpha > 1/2$ is quite similar (there are no poles inside the contour), hence it rests to consider $\alpha = \pm 1/2$. Again, both these cases are similar. Let $\alpha = 1/2$ and consider the contour composed by straight lines connecting the points $1-iT$, $1+iT$, $k+iT$, $k-iT$ with *k>1* and $T \to +\infty$ which semi-indents then point *z=1*. Then the repetition of the same that was done above gives $\frac{1}{\pi}\int_{-\infty}^{\infty}\cosh(bt/2)\operatorname{Re} F_\alpha(t)dt = 2 + 2\cos(b/2) - \operatorname{Re}((e^{-ib/2}+e^{ib/2})(1+2w(e^{ib})))$

and $\frac{1}{\pi}\int_{-\infty}^{\infty}\sinh(bt/2)\operatorname{Im} F_\alpha(t)dt = 2\sin(b/2) - \operatorname{Im}((-e^{-ib/2}+e^{ib/2})(1+2w(e^{ib})))$, where integrals are the principal values. Differentiating *2p* times and tending *b* to $\pi/2$ we get:

$$\frac{1}{\pi}\int_0^\infty t^{2p}\cosh(\pi t/4)\operatorname{Re} F_\alpha(t)dt = (-1)^p \cos(\pi/4)$$

$$\frac{1}{\pi}\int_0^\infty t^{2p}\sinh(\pi t/4)\operatorname{Im} F_\alpha(t)dt = (-1)^p \sin(\pi/4),$$ where the integrals do not contain any irregularities any more, and functions $t^p\cosh(\pi t/4)\operatorname{Re} F_\alpha(t)$, $t^p\sinh(\pi t/4)\operatorname{Im} F_\alpha(t)$ are continuous starting from *p=1* and *p=0* respectively. Thus Fejér theorem is again applicable and this finishes the proof.

I would like to finish the paper with the following remarks. First, it looks instructive to compare the situation which, in view of the proven Theorem, appears for the Riemann function zeroes with some known function. The trivial but, we believe, illuminative, example, is the "search" of zeroes for $\sinh(z^2)$ (or similar) function. Of course, $\sinh(x+iy)^2 = \sinh(x^2 - y^2)\cos(2xy) + i\cosh(x^2 - y^2)\sin(2xy)$ thus the real part of $\sinh(z^2)$ has zeroes along the straight lines $y = \pm x$ (they may be named *exceptional* lines similarly to exceptional line Re*s*=1/2 for the *F*-function where $\operatorname{Im} F = 0$) and along the hyperboles $y = \frac{\pi/4 + \pi n/2}{x}$, while its imaginary part has zeroes along the hyperboles $y = \frac{\pi n/2}{x}$. Of course, hyperboles $y = \frac{\pi/4 + \pi n/2}{x}$ and $y = \frac{\pi n/2}{x}$ have no intersections (we may say that some "special symmetry" is present, and for this simple example it is clear which one) and this is just the intersections of hyperboles $y = \frac{\pi n/2}{x}$ with the lines $y = \pm x$ which give zeroes of $\sinh(z^2)$ function; all of them lye on these two straight lines despite both real and imaginary parts of the function $\sinh(z^2)$ have infinitely many zeroes for any *x* (i.e. *Re*z) – again in an analogy with the function *F(z)*.

Second, it might be instructive to look for the proof of the theorems *how many* zeroes both real and imaginary parts of *F(z)* might have. It seems evident that for small enough $|\alpha|$ real part of *F* has at least *KT* zeroes, and apparently even *KTlnT* zeroes, on the segment $(1/2+\alpha, 1/2+\alpha+iT)$ (*K* is an appropriate constant) simply because for such a case $\operatorname{Re} F_\alpha(t) \cong |F_\alpha(t)|$. It seems that the proof of Hardy and Littlewood based on estimation of integrals $\int_{1/2+\alpha-iT}^{1/2+\alpha+iT} F(u) x^{-u/2} du$ over finite intervals [4] still holds with only minor modifications. Situation is more complicated with the imaginary part of *F* which for small $|\alpha|$ is small and the attempt to modify the proof meets serious difficulties. Could we believe that exactly the study of the solutions of an equation $\operatorname{Im} F_\alpha = 0$ for $\alpha \neq 0$, i.e. an equation $\int_1^\infty x^{-3/4}(x^{\alpha/2} - x^{-\alpha/2})\sin(\frac{t}{2}\ln x)w(x)dx = \frac{2\alpha t}{(\alpha^2 - 1/4 - t^2)^2 + 4\alpha^2 t^2}$, might shed some light on the Riemann hypothesis?